\documentclass[3p,times]{elsarticle}

\usepackage{ecrc}

\volume{}

\firstpage{1}

\journalname{Working Paper (Last Updated: June 25, 2014)}

\runauth{J. Chen}

\jnltitlelogo{Working Paper}

\CopyrightLine{Jeremy Chen 2014;}{Last Updated: June 25, 2014}

\usepackage{lineno}

\usepackage[figuresright]{rotating}

\usepackage{amsmath}
\usepackage{amssymb}
\usepackage{amsfonts}

\usepackage{enumerate}
\usepackage{tikz}
\usetikzlibrary{shapes}
\usetikzlibrary{arrows}
\usetikzlibrary{plotmarks}

\def\redTextAsRed{1}            	
\def\blueTextAsBlue{1}          	
\def\greenTextAsGreen{1}        	
\def\cyanTextAsCyan{1}			
\def\magentaTextAsMagenta{1}		

\def\commentsOn{1}				
\def\defaultCommentColour{magenta}  

\newcommand{\RedText}[1]{\ifthenelse{\redTextAsRed = 1}{{\color{red} #1}}{#1}}
\newcommand{\BlueText}[1]{\ifthenelse{\blueTextAsBlue = 1}{{\color{blue} #1}}{#1}}
\newcommand{\GreenText}[1]{\ifthenelse{\greenTextAsGreen = 1}{{\color{green} #1}}{#1}}
\newcommand{\MagentaText}[1]{\ifthenelse{\magentaTextAsMagenta = 1}{{\color{magenta} #1}}{#1}}
\newcommand{\CyanText}[1]{\ifthenelse{\cyanTextAsCyan = 1}{{\color{cyan} #1}}{#1}}

\newcommand{\CommentText}[1]{
    \ifthenelse{\commentsOn = 1}{{\color{\defaultCommentColour}{\noindent [{\bf Comment:} #1]}}}{}
}

\newcommand{\CommentTextN}[2]{
    \ifthenelse{\commentsOn = 1}{{\color{\defaultCommentColour}{\noindent [{\bf Comment (#1):} #2]}}}{}
}

\newcommand{\CommentTextC}[2]{
    \ifthenelse{\commentsOn = 1}{{\color{#1}{\noindent [{\bf Comment:} #2]}}}{}
}

\newcommand{\CommentTextCN}[3]{
    \ifthenelse{\commentsOn = 1}{{\color{#1}{\noindent [{\bf Comment (#2):} #3]}}}{}
}

\newcommand{\CommentTextLong}[1]{
    \ifthenelse{\commentsOn = 1}{{\color{\defaultCommentColour}{\noindent [{\bf --- Comment BEGIN ---}\newline\noindent #1 \ \newline\hspace*{\fill}{\bf --- Comment END ---}]}}}{}
}

\newcommand{\CommentTextLongN}[2]{
    \ifthenelse{\commentsOn = 1}{{\color{\defaultCommentColour}{\noindent [{\bf --- Comment (#1) BEGIN ---}\newline\noindent #2 \ \newline\hspace*{\fill}{\bf --- Comment (#1) END ---}]}}}{}
}

\newcommand{\CommentTextLongC}[2]{
    \ifthenelse{\commentsOn = 1}{{\color{#1}{\noindent [{\bf --- Comment BEGIN ---}\newline\noindent #2 \ \newline\hspace*{\fill}{\bf --- Comment END ---}]}}}{}
}

\newcommand{\CommentTextLongCN}[3]{
    \ifthenelse{\commentsOn = 1}{{\color{#1}{\noindent [{\bf --- Comment (#2) BEGIN ---}\newline\noindent #3 \ \newline\hspace*{\fill}{\bf --- Comment (#2) END ---}]}}}{}
}

\newcommand{\Ex}[1]{\mathbb{E}\left[#1\right]}

\newtheorem{theorem}{Theorem}[section]
\newtheorem{corollary}[theorem]{Corollary}
\newtheorem{lemma}[theorem]{Lemma}
\newtheorem{proposition}[theorem]{Proposition}

\usepackage{enumerate}
\usepackage{tikz}
\usetikzlibrary{shapes}
\usetikzlibrary{arrows}
\usetikzlibrary{plotmarks}

\begin{document}

\begin{frontmatter}



\dochead{}

\title{Weighing the ``Heaviest'' Polya Urn}


\author{Jeremy Chen}
\ead{jeremy.chen@nus.edu.sg}
\address{Department of Decision Sciences, National University of Singapore Business School\\15 Kent Ridge Drive, Singapore 119245}

\begin{abstract}

For the classical Polya urn model parametrized by exponent $\gamma$, the limit distributions of the fraction of balls in each bin are simple and well known when $\gamma<1$ (``egalitarian'') and when $\gamma>1$ (``winner takes all'').
In this note, we partially fill in the gap for $\gamma=1$, the critical point, by providing explicit analytical expressions for all the moments of the limit distribution of the fraction of balls in the bin with the most balls.

\end{abstract}


\end{frontmatter}

\linenumbers






\section{Preliminaries}\label{intro}

The Polya urn problem describes a well-studied family of random processes that have been fruitfully applied in diverse fields ranging from telecommunications to understanding self-organizing processes like network formation and herd behavior.
In the classical Polya urn problem, one begins with $d$ bins, each containing one ball. Additional balls arrive one at a time, and the probability that an arriving ball is placed in a given bin is proportional to $m^\gamma$, where $m$ is the number of balls in that bin.

In this note, we consider the case of $\gamma = 1$, which corresponds to a process of ``proportional preferential attachment'' and is a critical point with respect to the limit distribution of the fraction of balls in each bin. 
It is well known that for $\gamma < 1$ the fraction of balls in the ``heaviest'' bin (the bin with the most balls) tends to $1/d$, and for $\gamma > 1$ the fraction of balls in the ``heaviest'' bin tends to $1$. (See, for instance, surveys such as \citep{Chung2003} and \citep{Pemantle2007} or books such as \citep{JohnsonKotz1977} and \citep{Mahmoud2008}.) Unexpectedly, though scientists and engineers are interested in analogous quantities such as ``the size of the largest (biological) plague'', this question has not been explored for the case of $\gamma=1$.

The modest contribution of this note is to explore that very problem. We characterize the limit distribution of the fraction of balls in the ``heaviest'' bin for $\gamma=1$ by providing explicit analytical expressions for all its moments.


\section{The Main Result}\label{main}

Denote the number of balls in the ``heaviest'' urn (when there are $d$ urns), after a total of $n-d \ge 0$ balls are added, as $H_d(n)$. For integer $m \ge 0$ and integer $d \ge 1$, let
\begin{equation}\label{frac_def}
	M^{(m)}_d := \lim_{n\rightarrow\infty}\Ex{\left(\frac{H_d(n)}{n}\right)^m}.
\end{equation}
For notational convenience, let $M_0^{0} := 1$. We note that there are three equivalent ways of describing $M^{(m)}_d$:


\begin{proposition}[Limiting Moments of the Fraction of Balls in the ``Heaviest'' Urn]\label{prop_asymp_frac_exp}
	For $m\ge 1$ and $d \ge 2$,
	\begin{equation}\label{asymp_frac_moments_1}
		M^{(m)}_d
		= \sum_{k=0}^{m}{
			\frac{d-1}{d^k} 
			\frac{m!}{(m-k)!}
			\frac{(m+d-k-2)!}{(m+d-1)!}
			M^{(m-k)}_{d-1}
		},
	\end{equation}
	\begin{equation}\label{asymp_frac_moments_2}
		M^{(m)}_d
		= \frac{d-1}{m+d-1} M^{(m)}_{d-1} + \frac{m}{d(m+d-1)} M^{(m-1)}_d, {\rm\ and}
	\end{equation}
	\begin{equation}\label{asymp_frac_moments_3}
		M^{(m)}_d
		= \frac{m}{m+d-1} \sum_{j=1}^{d}{ \frac{(d-1)!}{j!} \frac{(m+j-2)!}{(m+d-2)!} M^{(m-1)}_j}.
	\end{equation}
\end{proposition}

This result is an explicit characterization in the sense that the computations to be done using either of equations (\ref{asymp_frac_moments_1}) thru (\ref{asymp_frac_moments_3}) are explicit (typically simple arithmetic) rather than implicit (e.g.: solving a system of equations). 
One readily observes that this is so because the next $M^{(m)}_d$ to be computed depends only on already computed moments. (It is clear that $M^{(0)}_d = 1$ for all $d\ge 1$ and $M^{(m)}_1 = 1$ for all $m\ge 1$.)

The first and second moments being typically of special interest, we may use equation (\ref{asymp_frac_moments_3}) to compute them:
\begin{corollary}[Limiting Mean/Second Moment of the Fraction of Balls in the ``Heaviest'' Urn]\label{cor_asymp_frac_mean_2M}
	\begin{equation}\label{asymp_frac_mean}
		\frac{1}{d}\log d \le M^{(1)}_d = \frac{1}{d}\sum_{k=1}^{d}{\frac{1}{k}} \le \frac{1}{d} (\log d + 1)
	\end{equation}
	and
	\begin{equation}\label{asymp_frac_2M}
		M^{(2)}_d = \frac{2}{d(d+1)} \sum_{j=1}^{d}{ M^{(1)}_j }
			= \frac{2}{d(d+1)} \sum_{j=1}^{d}{ \frac{1}{j} \sum_{k=1}^{j}\frac{1}{k} }.
	\end{equation}
\end{corollary}
Both ``scale like'' $d^{-1} \log d$.
Direct computation reveals that the coefficient of variation (the ratio of the standard deviation to the mean) of the fraction of balls in the ``heaviest'' urn rises to a peak of about 0.27 at $d=10$ after which it begins to decline with increading $d$ to a limiting value of $0$.

More generally, though the recurrences enable us to obtain expressions for any given moment by iterated substitution, there does not appear to be any special structure that enables drastic simplification. This is unfortunate.


Before proceeding to the prove Proposition \ref{prop_asymp_frac_exp}, we reproduce, for completeness, the following well known result on the distribution of the number of balls in each bin:

\begin{lemma}\label{lemma_pref_attach_uniform}
	The number of balls in the $d$ urns after $m$ additional balls are added is uniformly distributed on
\begin{equation}
	S_{(d,m)} := \{v \in \mathbb{Z}^d: v \ge e, v^T e = d+m\}.
\end{equation}
Furthermore,
\begin{equation}\label{eq_num_states}
	|S_{(d,m)}| = \frac{(m+d-1)!}{m!(d-1)!} = \left( m + d - 1 \atop d - 1\right).
\end{equation}
\end{lemma}
{\noindent\bf Proof of Lemma \ref{lemma_pref_attach_uniform}:}
Let the probability that $v_k$ balls are in the $k$-th urn (of $d$) after $m$ additional balls are added be $\pi_{(d,m)}(v)$. Clearly, $\pi_{(d,0)}(e) = 1$.
If $\pi_{(d,m)}(v) = \eta_{(d,m)}$ for some constant $\eta_{(d,m)}$ for all $v \in S_{(d,m)}$. Then, for all $v\in S_{(d,m+1)}$,
\begin{align}
	\pi_{(d,m+1)}(v)
		& = \sum_{v_k > 1} { \eta_{(d,m)} \frac{v_k - 1}{m+d} } \nonumber\\
		& = \sum_{k=1}^{d} { \eta_{(d,m)} \frac{v_k - 1}{m+d} } \nonumber\\
		& = \eta_{(d,m)} \frac{m+1}{m+d} \label{eq_num_states_inv} \\
		& =: \eta_{(d,m+1)} \nonumber
\end{align}
where the first equality follows from the dynamics of preferential attachment.
Now, $|S_{(d,m)}| = 1/\eta_{(d,m)}$, and $|S_{(d,0)}| = 1$. Therefore, by equation (\ref{eq_num_states_inv}),
\begin{equation*}
	S_{(d,m)} = 1\cdot\frac{d}{1}\cdot\frac{d+1}{2}\ldots \frac{m+d-1}{m}
\end{equation*}
and the proof is complete.
\hfill$\blacksquare$

Subsequently, Lemma \ref{lemma_pref_attach_uniform} and a simple partitioning of the set of possible outcomes ($S_{(d,m)}$, a discrete simplex) will be used to characterize the limiting distribution of the fraction of balls in the ``heaviest'' urn. We will also make use of the easily verifiable fact that:

\begin{lemma}\label{simple_lemma}
	For integer $a,b>0$ and real-valued $c>0$,
	\begin{equation}\label{simple_lemma_eq}
		\int_{0}^{c}{x^a (c-x)^b\;dx} = \frac{a! b!}{(a+b+1)!} c^{a+b+1}.
	\end{equation}
\end{lemma}


{\noindent\bf Proof of Proposition \ref{prop_asymp_frac_exp}:}
Clearly, $M^{(0)}_d = 1$ and $M^{(m)}_1 = 1$.
For cases where $m > 0$ and/or $d > 1$, an expression for the desired moment for finite $n$ will be constructed, and the limit as $n\rightarrow\infty$ evaluated.

Now, the set $S_{(d, (\alpha-1)d )}$ (for $\alpha\in\mathbb{N}$), as defined in Lemma \ref{lemma_pref_attach_uniform}, may be expressed as the following disjoint union:
\begin{equation*}
	S_{(d, (\alpha-1)d )} = \{\alpha e\} \cup 
		\bigcup_{\mu=1}^{\alpha - 1} \bigcup_{\tau=1}^{d - 1} T_{(\alpha, d, \mu, \tau)}
\end{equation*}
where
$
	T_{(\alpha, d, \mu, \tau)} := \{v \in \mathbb{Z}^d: v \ge \mu e, v^T e = \alpha d, \gamma(v, \mu) = \tau\}
$
and $\gamma(v,x) := | \{k : v_k = x\} |$ is the number of entries of the vector $v$ with the value $x$. This is so because there is a single vector in $S_{(d, (\alpha-1)d )}$ where all entries take the same value ($\alpha e$), and $\alpha-1$ other possible values of the smallest entry of vectors in $S_{(d, (\alpha-1)d )}$ (specifically, $1, 2, \ldots, \alpha - 1$). In each of the latter cases, the number of entries taking on the minimum value may range from $1$ to $d-1$.

Noting that vectors in $T_{(\alpha, d, \mu, \tau)}$ each have $\tau$ entries taking value $\mu$, and $d-\tau$ entries taking values strictly larger than $\mu$, the cardinality of $T_{(\alpha, d, \mu, \tau)}$, must be $|S_{(d-\tau, (\alpha-\mu)d - (d-\tau) )}|$ multiplied by the number of ways to pick the $\tau$ entries taking value $\mu$. Therefore, using Lemma \ref{lemma_pref_attach_uniform}, one may deduce that
\begin{equation}\label{cardinality_of_T}
	|T_{(\alpha, d, \mu, \tau)}| = |S_{(d-\tau, (\alpha-\mu)d - (d-\tau) )}| \left( d \atop \tau \right) = \left( (\alpha - \mu) d - 1 \atop d - \tau - 1\right) \left( d \atop \tau \right).
\end{equation}
Thus one obtains the identity
$
		\left( \alpha d - 1 \atop d - 1\right) = 1 + \sum_{\mu=1}^{\alpha - 1} \sum_{\tau=1}^{d - 1} \left( (\alpha - \mu) d - 1 \atop d - \tau - 1 \right) \left( d \atop \tau \right)
$
which, itself, may be proven directly by induction.

Now, equation (\ref{frac_def}) may be written equivalently as
\begin{equation}\label{asymp_frac_exp2}
	\Ex{H_d(n)^m} = M^{(m)}_d n^m + o(n^m)
\end{equation}
for integer $m\ge 1$.
Equation (\ref{asymp_frac_exp2}) clearly holds for $d=1$.
Now, suppose that equation (\ref{asymp_frac_exp2}) is true for $1, 2, \ldots, d-1$ urns. Using the fact that conditional on realizations being in $T_{(\alpha, d, \mu, \tau)}$, the uniform probability of outcomes implies that the distribution of $H_{d}(n)$ is identical to the (unconditional) distribution of $H_{d-\tau}((\alpha - \mu) d) + \mu$, the $m$-th moment of the number of balls in the ``heaviest'' urn is given by
\begin{align}
	\Ex{H_{d}(n)^m \big| T_{(\alpha, d, \mu, \tau)}}
	& = \Ex{(H_{d-\tau}((\alpha - \mu) d) + \mu)^m} \nonumber\\
	& = \sum_{k=0}^{m}{ \left( m \atop k \right)M^{(m-k)}_{d-\tau} ((\alpha - \mu) d)^{m-k} \mu^k} + o(\alpha^m) \label{eqn_conditional_moments}
\end{align}
following an application of equation (\ref{asymp_frac_exp2}).
Through equation (\ref{cardinality_of_T}), we obtain
\begin{equation}\label{the_expectation}
	\Ex{\left(\frac{H_{d}(\alpha d)}{\alpha d}\right)^m}
		= 
	\frac{1}{(\alpha d)^m}
	\frac{1}{\displaystyle\left( \alpha d - 1 \atop d - 1\right)}
	\left[
	\displaystyle\alpha^m + \sum_{\mu=1}^{\alpha - 1} \sum_{\tau=1}^{d - 1}
	\Ex{H_{d}(n)^m \big| T_{(\alpha, d, \mu, \tau)}}
	\left( (\alpha - \mu) d - 1 \atop d - \tau - 1 \right) \left( d \atop \tau \right)
	\right].
\end{equation}
Subsequently, taking limits,
\begin{align*}
	\lim_{\alpha\rightarrow\infty} \Ex{\left(\frac{H_{d}(\alpha d)}{\alpha d}\right)^m}
		& = \lim_{\alpha\rightarrow\infty} \frac{1}{(\alpha d)^m}
			\frac{1}{\displaystyle\left( \alpha d - 1 \atop d - 1\right)}
			\left[
			\displaystyle\alpha^m + \sum_{\mu=1}^{\alpha - 1} \sum_{\tau=1}^{d - 1}
			\Ex{H_{d}(n)^m \big| T_{(\alpha, d, \mu, \tau)}}
			\left( (\alpha - \mu) d - 1 \atop d - \tau - 1 \right)
			\left( d \atop \tau \right)
			\right]\\
		& = \lim_{\alpha\rightarrow\infty}
			\frac{
			\displaystyle\alpha^m + \sum_{\mu=1}^{\alpha - 1} \sum_{\tau=1}^{d - 1}
			\left[
				\sum_{k=0}^{m}{ \left( m \atop k \right) M^{(m-k)}_{d-\tau}((\alpha - \mu) d)^{m-k} \mu^k} + o(\alpha^m)
			\right]
			\left[
				\frac{ ((\alpha - \mu)d)^{d - \tau - 1} }{(d - \tau - 1)!} + o(\alpha^{d - \tau - 1})
			\right]
			\left( d \atop \tau \right)
			}
			{\displaystyle (\alpha d)^m \frac{(\alpha d)^{d-1}}{(d-1)!}} \\
		& = \lim_{\alpha\rightarrow\infty}
			\frac{
			\displaystyle \sum_{\mu=1}^{\alpha - 1} 
			\left(
				\left[
					\sum_{k=0}^{m}{ \left( m \atop k \right) M^{(m-k)}_{d-\tau}((\alpha - \mu) d)^{m-k} \mu^k} + o(\alpha^m)
				\right]
				\frac{ (\alpha - \mu)^{d-2} d^{d-2} }{(d-2)!} d
				+
				o(\alpha^{m+d-2})
			\right)
			}
			{\displaystyle (\alpha d)^m \frac{(\alpha d)^{d-1}}{(d-1)!}} \\
		& = \lim_{\alpha\rightarrow\infty}
			\frac{d-1}{d^m \alpha^{m+d-1}}
			\sum_{\mu=1}^{\alpha - 1}
			\left(
			\left[
				\sum_{k=0}^{m}{ \left( m \atop k \right) M^{(m-k)}_{d-1}((\alpha - \mu) d)^{m-k} \mu^k} + o(\alpha^m)
			\right] (\alpha - \mu)^{d-2}
			+
			o(\alpha^{m+d-2})
			\right)
			 \\
		& = \lim_{\alpha\rightarrow\infty}
			\frac{d-1}{d^m \alpha^{m+d-1}}
			\sum_{\mu=1}^{\alpha - 1} \left[
				\sum_{k=0}^{m}{ \left( m \atop k \right) M^{(m-k)}_{d-1}(\alpha - \mu)^{m+d-k-2} \mu^k d^{m-k}} + o(\alpha^{m+d-2})
			\right] \\
		& = \lim_{\alpha\rightarrow\infty}
			\frac{d-1}{d^m \alpha^{m+d-1}}
			\left[\int_{0}^{\alpha}{ 
				\sum_{k=0}^{m}{ \left( m \atop k \right) M^{(m-k)}_{d-1}(\alpha - \mu)^{m+d-k-2} \mu^k d^{m-k}}
			+ o(\alpha^{m+d-2}) \;d\mu } + o(\alpha^{m+d-1}) \right]\\
		& = \lim_{\alpha\rightarrow\infty}
			\frac{d-1}{d^m \alpha^{m+d-1}}
			\left[
				\sum_{k=0}^{m}{ \left( m \atop k \right) M^{(m-k)}_{d-1} \frac{(m+d-k-2)! k!}{(m+d-1)!} \alpha^{m+d-1} d^{m-k}}
			+ o(\alpha^{m+d-1})\right]\\
		& = \lim_{\alpha\rightarrow\infty}
			\left[
				\sum_{k=0}^{m}{\frac{d-1}{d^k}  \frac{m!}{(m-k)!} \frac{(m+d-k-2)!}{(m+d-1)!} M^{(m-k)}_{d-1}  }
			+ o(1)\right]\\
		& = 	\sum_{k=0}^{m}{\frac{d-1}{d^k}  \frac{m!}{(m-k)!} \frac{(m+d-k-2)!}{(m+d-1)!} M^{(m-k)}_{d-1}  }.
\end{align*}
The seventh equality follows from an application of Lemma \ref{simple_lemma} to evaluate the integral.

To complete the proof, we argue that the limit
$$
	\lim_{r\rightarrow\infty} \Ex{\left(\frac{H_{d}(n_r)}{n_r}\right)^m} =
\sum_{k=0}^{m}{\frac{d-1}{d^k}  \frac{m!}{(m-k)!} \frac{(m+d-k-2)!}{(m+d-1)!} M^{(m-k)}_{d-1}  }.
$$
also arises for all increasing positive integer sequences $\{n_r\}_{r=1}^{\infty}$
whose elements are greater than $d$ but are not all necessarily integral multiples of $d$.

Note that given any $n$, for $\beta_{x,d} := d \lfloor x/d\rfloor$, we have $\beta_{n,d} d \le n \le (\beta_{n,d}+1)d$. Since, in the right-hand-size of equation (\ref{the_expectation}), both the numerator and denominator are increasing in $\alpha$, one can construct upper and lower bounds by using either $\beta_{n,d}$ or $\beta_{n,d}+1$ accordingly in place of $\alpha$ in the numerator and denominator. For both bounds, the same leading order terms arise and hence the same limits. Therefore, with this sandwiching, we establish equation (\ref{asymp_frac_moments_1}).

With equation (\ref{asymp_frac_moments_1}), we may reduce the multiple term recurrence to a two term recurrence, yielding equation (\ref{asymp_frac_moments_2}), and then perform iterative substitution to obtain equation (\ref{asymp_frac_moments_3}). This proceeds as follows:
\begin{eqnarray*}
	M^{(m)}_d
	& = & \frac{d-1}{m+d-1} M^{(m)}_{d-1} + \frac{m}{d(m+d-1)} M^{(m-1)}_d \\
	& = & \prod_{j=2}^{d}{\frac{j-1}{m+j-1}} + \sum_{j=2}^{d-1}{\left(\prod_{k=j+1}^{d}\frac{k-1}{m+k-1}\right) \frac{m}{j(m+j-1)} M^{(m-1)}_j} + \frac{m}{d(m+d-1)} M^{(m-1)}_d \\
	& = & \frac{m!(d-1)!}{(m+d-1)!} + \sum_{j=2}^{d-1}{ m\frac{(d-1)!}{j!} \frac{(m+j-2)!}{(m+d-1)!} M^{(m-1)}_j} + \frac{m}{d(m+d-1)} M^{(m-1)}_d.
\end{eqnarray*}
(The second equality holds because $M^{(m-1)}_1=1$.) Simplifying yields equation (\ref{asymp_frac_moments_3}), as desired. 
\hfill$\blacksquare$

\section{Numerical Experiments
}\label{numerical_experiments}

Illustrations of Proposition \ref{prop_asymp_frac_exp} may be found in Figures \ref{fig_asymp_frac_exp_100} thru \ref{fig_asymp_frac_exp_20000} along with simulated average fractions of balls in the ``heaviest'' urn and quantiles obtained from simulation (for $n$ ranging from 100 to 20,000).
Those simulations suggest that the limiting values are good approximations even for systems with just thousands of balls added.



	\definecolor{orange}{rgb}{1,0.5,0}
	\definecolor{brown}{rgb}{0.75,0.25,0}

	\begin{figure}[H]\begin{center}
		\begin{tikzpicture}[y=4.05cm, x=.485cm]
			\draw (0,0) -- coordinate (x axis mid) (25,0);
		    	\draw (0,0) -- coordinate (y axis mid) (0,1);
		    	\foreach \x in {0,5,...,25} {
		     		\draw (\x,3pt) -- (\x,-3pt)
					node[anchor=north] {\x};
				}
		    	\foreach \y in {0,0.2,0.4,0.6,0.8,1} {
					\draw (3pt,\y) -- (-3pt,\y) 
		     			node[anchor=east] {\y}; 
				}
			\node[below=0.4cm] at (x axis mid) {$d$};
			\node[rotate=90, above=0.8cm] at (y axis mid) {$\displaystyle{H_d(n)} / n$};
		
			\draw[blue, thick] plot file {herding_asymp_exact_25.data};
			\draw[red, thick] plot[mark=*, only marks] file {herding_sim_100.data};
			\draw[brown, dashed] plot file {herding_sim_100_q20.data};
			\draw[brown, dashed] plot file {herding_sim_100_q80.data};
			\draw[orange, dashdotted] plot file {herding_sim_100_q05.data};
			\draw[orange, dashdotted] plot file {herding_sim_100_q95.data};

			\begin{scope}[shift={(9,0.985)}]
				\draw[blue, thick, yshift=0\baselineskip]
					plot[blue] (0,0) -- (1.5,0) 
					node[right]{Expectation: Proposition \ref{prop_asymp_frac_exp} ($n\rightarrow\infty$)};
				\draw[red, thick, yshift=-1\baselineskip]
					plot[mark=*, only marks] (0.75,0);
				\draw[red, thick, yshift=-1\baselineskip] (1.5,0) -- (1.5,0) 
					node[right]{Simulation Mean};
				\draw[brown, dashed, yshift=-2\baselineskip] (0,0) -- 
					plot[brown] (0,0) -- (1.5,0) 
					node[right]{Simulation Quantiles: 0.2 and 0.8};
				\draw[orange, dashdotted, yshift=-3\baselineskip] (0,0) -- 
					plot[orange] (0,0) -- (1.5,0) 
					node[right]{Simulation Quantiles: 0.05 and 0.95};
			\end{scope}
		\end{tikzpicture}
		\caption{Maximum Fraction of Balls in the ``Heaviest'' Urn: The Limiting Mean and Simulated Data ($n={\rm 100}$; 10,000 samples)}
		\label{fig_asymp_frac_exp_100}
	\end{center}\end{figure}
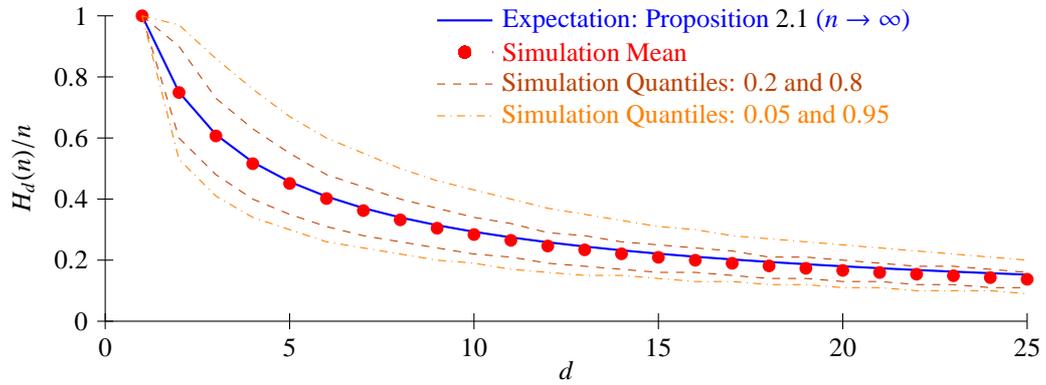

	\begin{figure}[H]\begin{center}
		\begin{tikzpicture}[y=4.05cm, x=.485cm]
			\draw (0,0) -- coordinate (x axis mid) (25,0);
		    	\draw (0,0) -- coordinate (y axis mid) (0,1);
		    	\foreach \x in {0,5,...,25} {
		     		\draw (\x,3pt) -- (\x,-3pt)
					node[anchor=north] {\x};
				}
		    	\foreach \y in {0,0.2,0.4,0.6,0.8,1} {
					\draw (3pt,\y) -- (-3pt,\y) 
		     			node[anchor=east] {\y}; 
				}
			\node[below=0.4cm] at (x axis mid) {$d$};
			\node[rotate=90, above=0.8cm] at (y axis mid) {$\displaystyle{H_d(n)} / n$};
		
			\draw[blue, thick] plot file {herding_asymp_exact_25.data};
			\draw[red, thick] plot[mark=*, only marks] file {herding_sim_250.data};
			\draw[brown, dashed] plot file {herding_sim_250_q20.data};
			\draw[brown, dashed] plot file {herding_sim_250_q80.data};
			\draw[orange, dashdotted] plot file {herding_sim_250_q05.data};
			\draw[orange, dashdotted] plot file {herding_sim_250_q95.data};

			\begin{scope}[shift={(9,0.985)}]
				\draw[blue, thick, yshift=0\baselineskip]
					plot[blue] (0,0) -- (1.5,0) 
					node[right]{Expectation: Proposition \ref{prop_asymp_frac_exp} ($n\rightarrow\infty$)};
				\draw[red, thick, yshift=-1\baselineskip]
					plot[mark=*, only marks] (0.75,0);
				\draw[red, thick, yshift=-1\baselineskip] (1.5,0) -- (1.5,0) 
					node[right]{Simulation Mean};
				\draw[brown, dashed, yshift=-2\baselineskip] (0,0) -- 
					plot[brown] (0,0) -- (1.5,0) 
					node[right]{Simulation Quantiles: 0.2 and 0.8};
				\draw[orange, dashdotted, yshift=-3\baselineskip] (0,0) -- 
					plot[orange] (0,0) -- (1.5,0) 
					node[right]{Simulation Quantiles: 0.05 and 0.95};
			\end{scope}
		\end{tikzpicture}
		\caption{Maximum Fraction of Balls in the ``Heaviest'' Urn: The Limiting Mean and Simulated Data ($n={\rm 250}$; 10,000 samples)}
		\label{fig_asymp_frac_exp_250}
	\end{center}\end{figure}
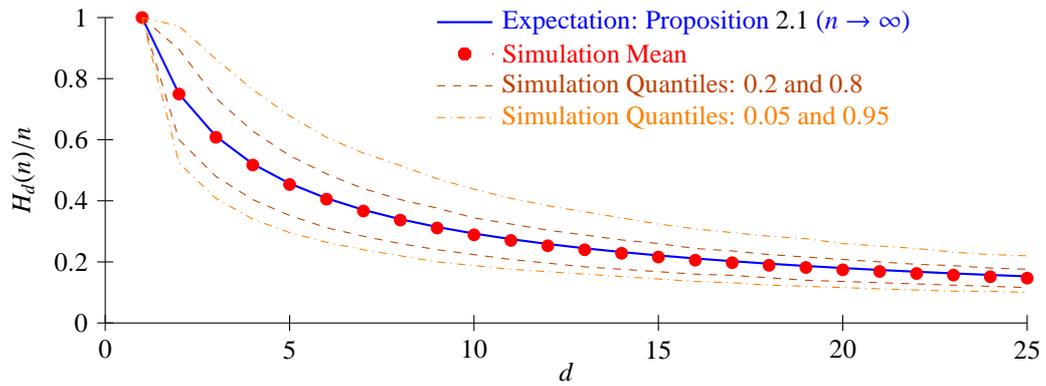

	\begin{figure}[H]\begin{center}
		\begin{tikzpicture}[y=4.05cm, x=.485cm]
			\draw (0,0) -- coordinate (x axis mid) (25,0);
		    	\draw (0,0) -- coordinate (y axis mid) (0,1);
		    	\foreach \x in {0,5,...,25} {
		     		\draw (\x,3pt) -- (\x,-3pt)
					node[anchor=north] {\x};
				}
		    	\foreach \y in {0,0.2,0.4,0.6,0.8,1} {
					\draw (3pt,\y) -- (-3pt,\y) 
		     			node[anchor=east] {\y}; 
				}
			\node[below=0.4cm] at (x axis mid) {$d$};
			\node[rotate=90, above=0.8cm] at (y axis mid) {$\displaystyle{H_d(n)} / n$};
		
			\draw[blue, thick] plot file {herding_asymp_exact_25.data};
			\draw[red, thick] plot[mark=*, only marks] file {herding_sim_500.data};
			\draw[brown, dashed] plot file {herding_sim_500_q20.data};
			\draw[brown, dashed] plot file {herding_sim_500_q80.data};
			\draw[orange, dashdotted] plot file {herding_sim_500_q05.data};
			\draw[orange, dashdotted] plot file {herding_sim_500_q95.data};

			\begin{scope}[shift={(9,0.985)}]
				\draw[blue, thick, yshift=0\baselineskip]
					plot[blue] (0,0) -- (1.5,0) 
					node[right]{Expectation: Proposition \ref{prop_asymp_frac_exp} ($n\rightarrow\infty$)};
				\draw[red, thick, yshift=-1\baselineskip]
					plot[mark=*, only marks] (0.75,0);
				\draw[red, thick, yshift=-1\baselineskip] (1.5,0) -- (1.5,0) 
					node[right]{Simulation Mean};
				\draw[brown, dashed, yshift=-2\baselineskip] (0,0) -- 
					plot[brown] (0,0) -- (1.5,0) 
					node[right]{Simulation Quantiles: 0.2 and 0.8};
				\draw[orange, dashdotted, yshift=-3\baselineskip] (0,0) -- 
					plot[orange] (0,0) -- (1.5,0) 
					node[right]{Simulation Quantiles: 0.05 and 0.95};
			\end{scope}
		\end{tikzpicture}
		\caption{Maximum Fraction of Balls in the ``Heaviest'' Urn: The Limiting Mean and Simulated Data ($n={\rm 500}$; 10,000 samples)}
		\label{fig_asymp_frac_exp_500}
	\end{center}\end{figure}
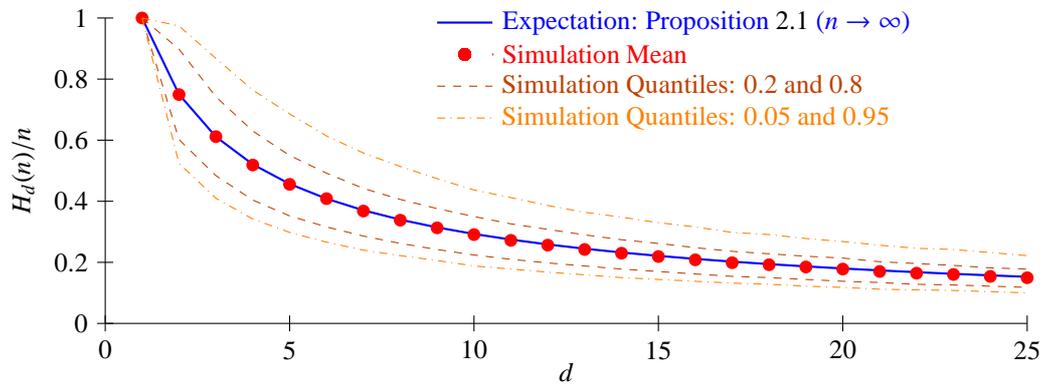

	\begin{figure}[H]\begin{center}
		\begin{tikzpicture}[y=4.05cm, x=.485cm]
			\draw (0,0) -- coordinate (x axis mid) (25,0);
		    	\draw (0,0) -- coordinate (y axis mid) (0,1);
		    	\foreach \x in {0,5,...,25} {
		     		\draw (\x,3pt) -- (\x,-3pt)
					node[anchor=north] {\x};
				}
		    	\foreach \y in {0,0.2,0.4,0.6,0.8,1} {
					\draw (3pt,\y) -- (-3pt,\y) 
		     			node[anchor=east] {\y}; 
				}
			\node[below=0.4cm] at (x axis mid) {$d$};
			\node[rotate=90, above=0.8cm] at (y axis mid) {$\displaystyle{H_d(n)} / n$};
		
			\draw[blue, thick] plot file {herding_asymp_exact_25.data};
			\draw[red, thick] plot[mark=*, only marks] file {herding_sim_1000.data};
			\draw[brown, dashed] plot file {herding_sim_1000_q20.data};
			\draw[brown, dashed] plot file {herding_sim_1000_q80.data};
			\draw[orange, dashdotted] plot file {herding_sim_1000_q05.data};
			\draw[orange, dashdotted] plot file {herding_sim_1000_q95.data};

			\begin{scope}[shift={(9,0.985)}]
				\draw[blue, thick, yshift=0\baselineskip]
					plot[blue] (0,0) -- (1.5,0) 
					node[right]{Expectation: Proposition \ref{prop_asymp_frac_exp} ($n\rightarrow\infty$)};
				\draw[red, thick, yshift=-1\baselineskip]
					plot[mark=*, only marks] (0.75,0);
				\draw[red, thick, yshift=-1\baselineskip] (1.5,0) -- (1.5,0) 
					node[right]{Simulation Mean};
				\draw[brown, dashed, yshift=-2\baselineskip] (0,0) -- 
					plot[brown] (0,0) -- (1.5,0) 
					node[right]{Simulation Quantiles: 0.2 and 0.8};
				\draw[orange, dashdotted, yshift=-3\baselineskip] (0,0) -- 
					plot[orange] (0,0) -- (1.5,0) 
					node[right]{Simulation Quantiles: 0.05 and 0.95};
			\end{scope}
		\end{tikzpicture}
		\caption{Maximum Fraction of Balls in the ``Heaviest'' Urn: The Limiting Mean and Simulated Data ($n={\rm 1,000}$; 10,000 samples)}
		\label{fig_asymp_frac_exp_1000}
	\end{center}\end{figure}
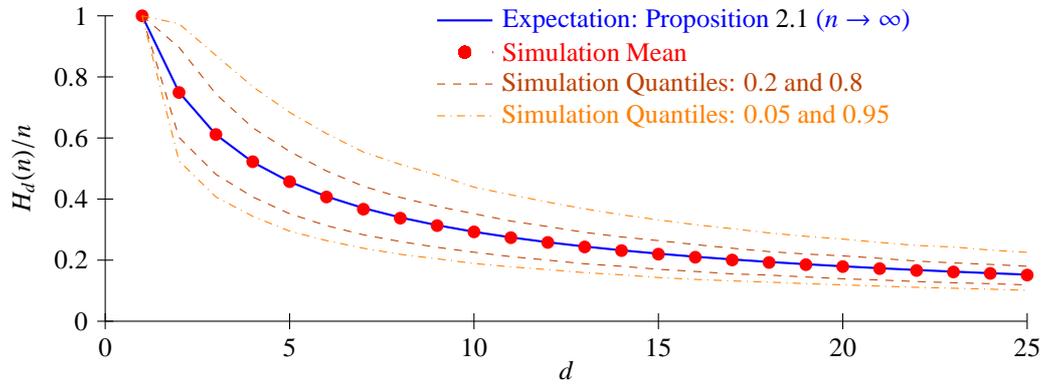

	\begin{figure}[H]\begin{center}
		\begin{tikzpicture}[y=4.05cm, x=.485cm]
			\draw (0,0) -- coordinate (x axis mid) (25,0);
		    	\draw (0,0) -- coordinate (y axis mid) (0,1);
		    	\foreach \x in {0,5,...,25} {
		     		\draw (\x,3pt) -- (\x,-3pt)
					node[anchor=north] {\x};
				}
		    	\foreach \y in {0,0.2,0.4,0.6,0.8,1} {
					\draw (3pt,\y) -- (-3pt,\y) 
		     			node[anchor=east] {\y}; 
				}
			\node[below=0.4cm] at (x axis mid) {$d$};
			\node[rotate=90, above=0.8cm] at (y axis mid) {$\displaystyle{H_d(n)} / n$};
		
			\draw[blue, thick] plot file {herding_asymp_exact_25.data};
			\draw[red, thick] plot[mark=*, only marks] file {herding_sim_2000.data};
			\draw[brown, dashed] plot file {herding_sim_2000_q20.data};
			\draw[brown, dashed] plot file {herding_sim_2000_q80.data};
			\draw[orange, dashdotted] plot file {herding_sim_2000_q05.data};
			\draw[orange, dashdotted] plot file {herding_sim_2000_q95.data};

			\begin{scope}[shift={(9,0.985)}]
				\draw[blue, thick, yshift=0\baselineskip]
					plot[blue] (0,0) -- (1.5,0) 
					node[right]{Expectation: Proposition \ref{prop_asymp_frac_exp} ($n\rightarrow\infty$)};
				\draw[red, thick, yshift=-1\baselineskip]
					plot[mark=*, only marks] (0.75,0);
				\draw[red, thick, yshift=-1\baselineskip] (1.5,0) -- (1.5,0) 
					node[right]{Simulation Mean};
				\draw[brown, dashed, yshift=-2\baselineskip] (0,0) -- 
					plot[brown] (0,0) -- (1.5,0) 
					node[right]{Simulation Quantiles: 0.2 and 0.8};
				\draw[orange, dashdotted, yshift=-3\baselineskip] (0,0) -- 
					plot[orange] (0,0) -- (1.5,0) 
					node[right]{Simulation Quantiles: 0.05 and 0.95};
			\end{scope}
		\end{tikzpicture}
		\caption{Maximum Fraction of Balls in the ``Heaviest'' Urn: The Limiting Mean and Simulated Data ($n={\rm 2,000}$; 10,000 samples)}
		\label{fig_asymp_frac_exp_2000}
	\end{center}\end{figure}
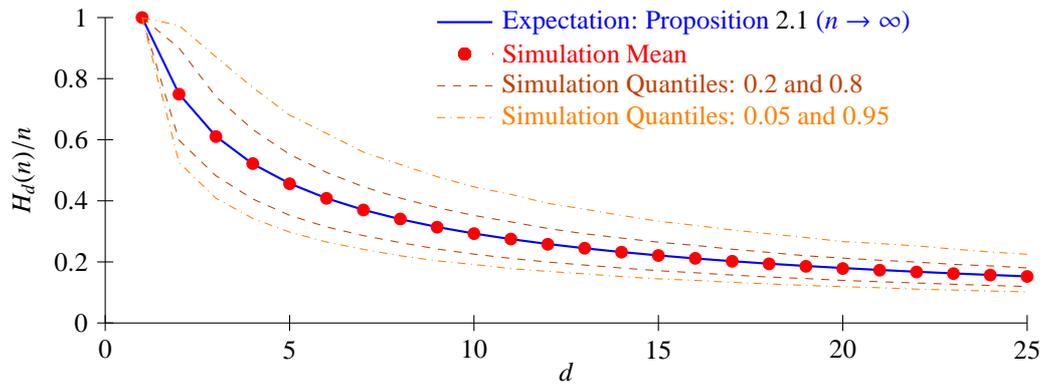

	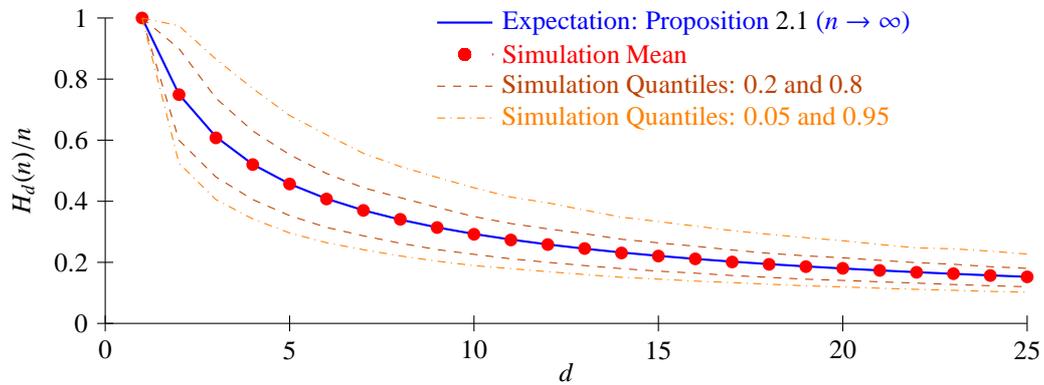
\begin{figure}[H]\begin{center}
		\begin{tikzpicture}[y=4.05cm, x=.485cm]
			\draw (0,0) -- coordinate (x axis mid) (25,0);
		    	\draw (0,0) -- coordinate (y axis mid) (0,1);
		    	\foreach \x in {0,5,...,25} {
		     		\draw (\x,3pt) -- (\x,-3pt)
					node[anchor=north] {\x};
				}
		    	\foreach \y in {0,0.2,0.4,0.6,0.8,1} {
					\draw (3pt,\y) -- (-3pt,\y) 
		     			node[anchor=east] {\y}; 
				}
			\node[below=0.4cm] at (x axis mid) {$d$};
			\node[rotate=90, above=0.8cm] at (y axis mid) {$\displaystyle{H_d(n)} / n$};
		
			\draw[blue, thick] plot file {herding_asymp_exact_25.data};
			\draw[red, thick] plot[mark=*, only marks] file {herding_sim_5000.data};
			\draw[brown, dashed] plot file {herding_sim_5000_q20.data};
			\draw[brown, dashed] plot file {herding_sim_5000_q80.data};
			\draw[orange, dashdotted] plot file {herding_sim_5000_q05.data};
			\draw[orange, dashdotted] plot file {herding_sim_5000_q95.data};

			\begin{scope}[shift={(9,0.985)}]
				\draw[blue, thick, yshift=0\baselineskip]
					plot[blue] (0,0) -- (1.5,0) 
					node[right]{Expectation: Proposition \ref{prop_asymp_frac_exp} ($n\rightarrow\infty$)};
				\draw[red, thick, yshift=-1\baselineskip]
					plot[mark=*, only marks] (0.75,0);
				\draw[red, thick, yshift=-1\baselineskip] (1.5,0) -- (1.5,0) 
					node[right]{Simulation Mean};
				\draw[brown, dashed, yshift=-2\baselineskip] (0,0) -- 
					plot[brown] (0,0) -- (1.5,0) 
					node[right]{Simulation Quantiles: 0.2 and 0.8};
				\draw[orange, dashdotted, yshift=-3\baselineskip] (0,0) -- 
					plot[orange] (0,0) -- (1.5,0) 
					node[right]{Simulation Quantiles: 0.05 and 0.95};
			\end{scope}
		\end{tikzpicture}
		\caption{Maximum Fraction of Balls in the ``Heaviest'' Urn: The Limiting Mean and Simulated Data ($n={\rm 5,000}$; 10,000 samples)}
		\label{fig_asymp_frac_exp_2000}
	\end{center}\end{figure}

	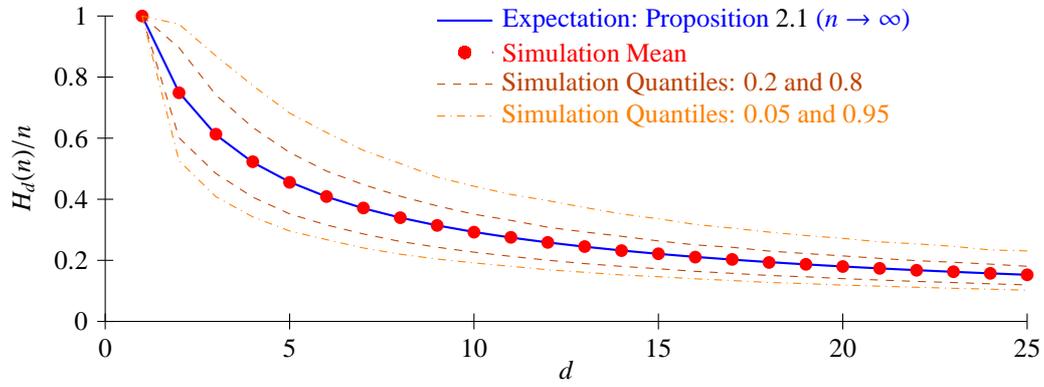
\begin{figure}[H]\begin{center}
		\begin{tikzpicture}[y=4.05cm, x=.485cm]
			\draw (0,0) -- coordinate (x axis mid) (25,0);
		    	\draw (0,0) -- coordinate (y axis mid) (0,1);
		    	\foreach \x in {0,5,...,25} {
		     		\draw (\x,3pt) -- (\x,-3pt)
					node[anchor=north] {\x};
				}
		    	\foreach \y in {0,0.2,0.4,0.6,0.8,1} {
					\draw (3pt,\y) -- (-3pt,\y) 
		     			node[anchor=east] {\y}; 
				}
			\node[below=0.4cm] at (x axis mid) {$d$};
			\node[rotate=90, above=0.8cm] at (y axis mid) {$\displaystyle{H_d(n)} / n$};
		
			\draw[blue, thick] plot file {herding_asymp_exact_25.data};
			\draw[red, thick] plot[mark=*, only marks] file {herding_sim_10000.data};
			\draw[brown, dashed] plot file {herding_sim_10000_q20.data};
			\draw[brown, dashed] plot file {herding_sim_10000_q80.data};
			\draw[orange, dashdotted] plot file {herding_sim_10000_q05.data};
			\draw[orange, dashdotted] plot file {herding_sim_10000_q95.data};

			\begin{scope}[shift={(9,0.985)}]
				\draw[blue, thick, yshift=0\baselineskip]
					plot[blue] (0,0) -- (1.5,0) 
					node[right]{Expectation: Proposition \ref{prop_asymp_frac_exp} ($n\rightarrow\infty$)};
				\draw[red, thick, yshift=-1\baselineskip]
					plot[mark=*, only marks] (0.75,0);
				\draw[red, thick, yshift=-1\baselineskip] (1.5,0) -- (1.5,0) 
					node[right]{Simulation Mean};
				\draw[brown, dashed, yshift=-2\baselineskip] (0,0) -- 
					plot[brown] (0,0) -- (1.5,0) 
					node[right]{Simulation Quantiles: 0.2 and 0.8};
				\draw[orange, dashdotted, yshift=-3\baselineskip] (0,0) -- 
					plot[orange] (0,0) -- (1.5,0) 
					node[right]{Simulation Quantiles: 0.05 and 0.95};
			\end{scope}
		\end{tikzpicture}
		\caption{Maximum Fraction of Balls in the ``Heaviest'' Urn: The Limiting Mean and Simulated Data ($n={\rm 10,000}$; 10,000 samples)}
		\label{fig_asymp_frac_exp_10000}
	\end{center}\end{figure}

	\begin{figure}[H]\begin{center}
		\begin{tikzpicture}[y=4.05cm, x=.485cm]
			\draw (0,0) -- coordinate (x axis mid) (25,0);
		    	\draw (0,0) -- coordinate (y axis mid) (0,1);
		    	\foreach \x in {0,5,...,25} {
		     		\draw (\x,3pt) -- (\x,-3pt)
					node[anchor=north] {\x};
				}
		    	\foreach \y in {0,0.2,0.4,0.6,0.8,1} {
					\draw (3pt,\y) -- (-3pt,\y) 
		     			node[anchor=east] {\y}; 
				}
			\node[below=0.4cm] at (x axis mid) {$d$};
			\node[rotate=90, above=0.8cm] at (y axis mid) {$\displaystyle{H_d(n)} / n$};
		
			\draw[blue, thick] plot file {herding_asymp_exact_25.data};
			\draw[red, thick] plot[mark=*, only marks] file {herding_sim_20000.data};
			\draw[brown, dashed] plot file {herding_sim_20000_q20.data};
			\draw[brown, dashed] plot file {herding_sim_20000_q80.data};
			\draw[orange, dashdotted] plot file {herding_sim_20000_q05.data};
			\draw[orange, dashdotted] plot file {herding_sim_20000_q95.data};

			\begin{scope}[shift={(9,0.985)}]
				\draw[blue, thick, yshift=0\baselineskip]
					plot[blue] (0,0) -- (1.5,0) 
					node[right]{Expectation: Proposition \ref{prop_asymp_frac_exp} ($n\rightarrow\infty$)};
				\draw[red, thick, yshift=-1\baselineskip]
					plot[mark=*, only marks] (0.75,0);
				\draw[red, thick, yshift=-1\baselineskip] (1.5,0) -- (1.5,0) 
					node[right]{Simulation Mean};
				\draw[brown, dashed, yshift=-2\baselineskip] (0,0) -- 
					plot[brown] (0,0) -- (1.5,0) 
					node[right]{Simulation Quantiles: 0.2 and 0.8};
				\draw[orange, dashdotted, yshift=-3\baselineskip] (0,0) -- 
					plot[orange] (0,0) -- (1.5,0) 
					node[right]{Simulation Quantiles: 0.05 and 0.95};
			\end{scope}
		\end{tikzpicture}
		\caption{Maximum Fraction of Balls in the ``Heaviest'' Urn: The Limiting Mean and Simulated Data ($n={\rm 20,000}$; 10,000 samples)}
		\label{fig_asymp_frac_exp_20000}
	\end{center}\end{figure}
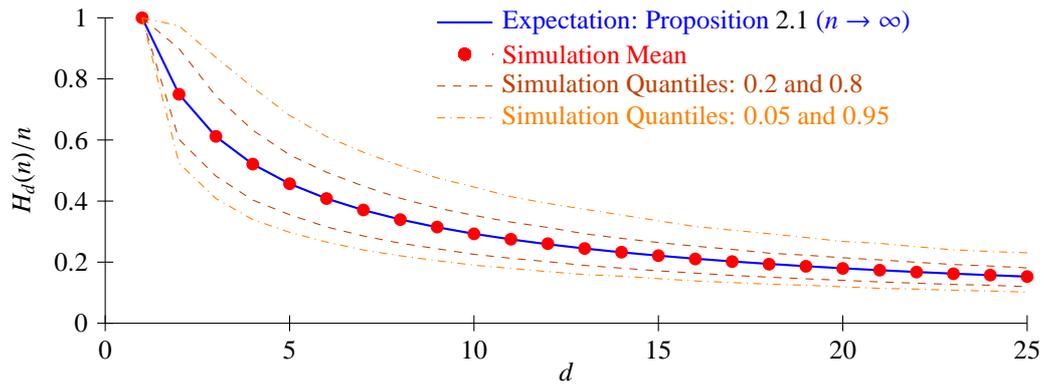


\section{Herding by Preferential Attachment}

For the author, the particular question of what fraction of the balls, in an instance of a Polya urn problem, would end up in the urn with the most balls arose from studying a discrete choice model of connected agents whose choices affect those of others \citep{SNatCPM2013WP}.
In that model, an agent would probabilistically select a choice based on his/her own preferences or adopt the choice of some other agent (i.e.: in ``steady state'', the former agent makes choices in accordance with the choice probabilities of the other agent).
By imposing specific decision making dynamics on that ``steady state model'', a stochastic model was obtained whereby, in any sample path, the realized choice of any given agent might be traced ``up the tree of choice adoption'' to some ``decisive agent'' (an agent making his/her choices based on his/her own preferences rather than by adopting that of another agent). As such, each outcome would correspond to a ``forest'' of ``trees'' each containing exactly one ``decisive agent'' and rooted at that agent.
Let us denote all the agents connected by choice adoption decisions to each ``decisive agent'' to be the ``herd'' associated with that ``decisive agent''.

To gain some insight on herd behavior, the model might be simplified to one where, conditional on the ``decisive agents'', any agent adopting the choice of another would select uniformly among the other agents. Further conditioning on there being $d$ ``decisive agents'', one obtains an instance of the Polya urn model with $d$ urns.
Using the largest herd as a proxy for the extent of herding behavior, Proposition \ref{prop_asymp_frac_exp} (and Corollary \ref{cor_asymp_frac_mean_2M}, in particular) inform one about the relationship between herding and the ``decisiveness'' of agents. Furthermore, should one desire to evaluate the ``risk'' associated with herding as the expectation of some (smooth) function of the fraction of agents in the largest herd, given that the moments are decreasing and bounded between $0$ and $1$ and that numerical experiments indicate the limiting distributions are ``close'' for reasonably sized populations of agents, good Taylor-expansion-based approximations can be computed.






\bibliographystyle{elsarticle-num}
\bibliography{bib}







\end{document}